\documentclass[a4paper,11pt]{amsart} 
\usepackage{amssymb,amsfonts,amsxtra}
\usepackage[all]{xy} \usepackage{fullpage}
\newtheorem{theorem}{Theorem}[section]
\newtheorem{cor}[theorem]{Corollary}

\newtheorem{prop}[theorem]{Proposition}

\newtheorem{example}[theorem]{Example}

\DeclareMathOperator{\Der}{Der}
\DeclareMathOperator{\ad}{ad}

\newcommand{\noproof}{\begin{flushright} \ensuremath{\square} \end{flushright}} \newtheorem{defi}[theorem]{Definition}
\newtheorem{rem}[theorem]{Remark} 
\numberwithin{equation}{section}
\subjclass[2000]{55P62, 13D03, 57T30}
\thanks{This research was partially supported by the EPSRC grant
No. GR/SO7148/01} \begin{document} \title[Hochschild cohomology
and moduli spaces...] {The Stasheff model of a simply-connected manifold and the string bracket} \author{A.Lazarev}\address{Mathematics Department, University of Bristol, Bristol,
BS8 1TW, England.} \email{a.lazarev@bristol.ac.uk}
\keywords{Rational homotopy, $C_\infty$-algebra, symplectic structure, string topology}
\begin{abstract}
We revisit Stasheff's construction of a minimal Lie-Quillen model of a simply-connected closed manifold $M$ using the language of infinity-algebras.
 This model is then used to construct a graded Lie bracket on the equivariant homology of the free loop space of $M$ minus a point similar to the Chas-Sullivan string bracket. 
\end{abstract}

\maketitle 

\section{Introduction}
In this paper we construct a string-type bracket on the $S^1$-equivariant reduced homology of the loop space of a simply-connected closed manifold $M$ with a puncture. In fact, we only need $M$ to be a rational Poincar\'e duality space of dimension $n$;  removing a point corresponds to passing to the $n-1$-skeleton of $M$. It seems likely that our construction is compatible with that of Sullivan and Chas \cite{SC} under the inclusion $\dot{M}:=M\setminus point \hookrightarrow M$, however this issue is not considered here. 

Our main tool is the notion of a symplectic infinity-algebra, a.k.a. infinity-algebra with an invariant inner product introduced by Kontsevich \cite{kontfd} and studied in detail in \cite{HL}. This is simply a homotopy invariant version of a (graded) Frobenius algebra. Since cohomology rings of Poincar\'e duality spaces are graded Frobenius algebras the appearance of symplectic infinity-algebras is not unexpected. Note that in the simply-connected case the approach using infinity algebras is a more or less tautological, albeit convenient, reformulation of the more traditional one via the Sullivan and Lie-Quillen models. 

 Consider a minimal Lie-Quillen model (or Quillen model for short) of a simply-connected manifold $M$; recall that it has the \emph{reduced} homology of $M$ as its underlying space, thus it does not support  the Poincar\'e duality form. To restore it one could either add a unit which results in what we call a \emph{contractible Quillen model} of $M$, or to remove the top class which corresponds to making a puncture in $M$.   Using the results of \cite{HL} we show that there exists a contractible Quillen model of $M$ that is a symplectic $C_\infty$-algebra. Similarly there is a Quillen model for $\dot{M}$ that is a symplectic $C_\infty$-algebra. The latter was essentially  constructed by Stasheff in \cite{Sta} and so we call it a \emph{Stasheff model}.

To construct the string bracket on the equivariant homology of the loops on $\dot{M}$ we use the connection of this homology with cyclic cohomology of the cochain algebra of $\dot{M}$ cf. \cite{jones}.

The paper is organized as follows. In section 2 we recall the definitions and basic facts about infinity algebras following  \cite{GJ}, \cite{Laz} and especially \cite{HL} and relate them to rational homotopy theory. In section 3 we consider cyclic cohomology of infinity-algebras. In section 4 we introduce symplectic infinity-algebras and use them to construct models for simply-connected Poincar\'e duality spaces. The string bracket is constructed in section 5.        
\subsection{Notation and conventions} We work over a fixed field $k$ of characteristic zero; all homology and cohomology groups are taken with coefficients in $k$. Whenever we talk about differential graded models of topological spaces $k$ is understood to be the field of rational numbers. The terms `differential graded algebra' and `differential graded Lie algebra' will be abbreviated as `dga' and `dgla' respectively. The $k$-dual to a graded vector space $V$ will be denoted by $V^*$ whilst the (homological or cohomological) grading will be indicated by an upper or lower bullet$~\bullet$. We will denote by $TV$ and $LV$ respectively the tensor algebra and the free Lie algebra on a graded vector space $V$. Their completions will be denoted by $\hat{T}V$ and $\hat{L}V$ respectively. The spaces of noncommutative power series or Lie series in indeterminates $x_i$ will be denoted by $k\langle\langle x_1,x_2,\ldots\rangle\rangle=k\langle\langle \bf x\rangle\rangle$ and by $k\{\{ x_1,x_2,\ldots \}\}=k\{\{ \bf x\}\}$ respectively. The suspension of a graded vector space $V^\bullet$ is defined as $\Sigma V^\bullet:=V^{\bullet+1}$. For a graded space or an algebra $V$ supplied with an augmentation $V\rightarrow k$ we denote by $V_+$ the kernel of the augmentation.

The $n-1$-skeleton of a rational Poincar\'e duality space $M$ of dimension $n$ will be denoted by $\dot{M}$.
\subsection{Acknowledgement} The author wishes to express his appreciation to M. Aubry, J.-L. Lemaire and J. Stasheff for useful discussions.
 
\section{Infinity-algebras and rational homotopy theory}
 \subsection{Generalities on infinity-algebras.} Recall that an $A_\infty$-algebra structure on a graded vector space $V$ is a continuous derivation $m$ of the completed tensor algebra $\hat{T}\Sigma V^*$ of homological degree $-1$, having square zero and vanishing at zero. Let us choose a basis $\{t_i\}$ in $\Sigma V^*$. Any element in $T\Sigma V^*$ is a noncommutative power series $f(t_1,t_2,\ldots)$ in the indeterminates $t_i$. 
Then we could write $m$ as $m=m_1+m_2+\ldots$ where $m_i=\sum_l f^l_i({\bf t})\partial_{t_l}$ where $f^l_i$ is a (possibly infinite) linear combination of monomials having wordlenth $l$. The condition $m^2=0$ implies that $m_1^2=0$ and so $m_1$ determines a differential on $V$. If $m_1=0$, i.e. if the differential $m$ is decomposable, we say that the $A_\infty$-structure $m$ is \emph{minimal}. In this case the quadratic part $m_2$ of $m$ determines an associative multiplication on $V$.

Next, a (continuous) derivation of $\hat{L}(\Sigma V^*)$   is a $C_\infty$-structure on $V$. It is clear that a $C_\infty$-algebra is a special case of an $A_\infty$-algebra.  The quadratic term $m_2$ of a minimal $C_\infty$-algebra  determines a \emph{commutative} product on $V$.          

Given two $A_\infty$-structures ($C_\infty$-structures) $m^V$ and $m^U$ on $V$ and $U$ an $A_\infty$-morphism ($C_\infty$-morphism) $V\rightarrow U$ is a continuous algebra homomorphism $\phi:\hat{T}\Sigma U^*\rightarrow \hat{T}\Sigma V^*$ ($\phi:\hat{L} \Sigma U^*\rightarrow \hat{L}\Sigma V^*$) for which  $\phi\circ m^U=m^U\circ V$. Such a  map $\phi$ could always be written as $\phi=\phi_1+\phi_2+\ldots$ where $\phi_i$ is a morphism raising the wordlength (or bracket length) by $i-1$. In particular, $\phi_1$ could be thought of as a linear map $\Sigma U^*\rightarrow \Sigma V^*$. We say that $\phi$ is a weak equivalence if $\phi_1$ determines a quasi-isomorphism between $\Sigma U^*$ and $\Sigma V^*$ considered as complexes with differentials $m^U_1$ and $m^V_1$. A weak equivalence between minimal infinity-algebras is always an isomorphism.

Given a differential graded algebra $V$ its cobar-construction $T\Sigma V^*$ could be considered as an $A_\infty$-algebra of a special sort. Indeed, the differential $m$ on $T\Sigma V^*$ is a sum of $m_1$ and $m_2$ which correspond to the differential and product on $V$ respectively. If $V$ is commutative then $T\Sigma V^*$ is in fact a differential graded Hopf algebra which gives rise to a $C_\infty$-algebra after taking the primitives. Kadeishvili's theorem states that any $A_\infty$-algebra (in particular the one corresponding to a cobar-construction of a differential graded algebra) admits a minimal model, i.e. a minimal $A_\infty$-algebra weakly equivalent to it. The analogue of this theorem is also valid in the $C_\infty$-case cf. \cite{HL}. 
\subsection{Adjoining a unit to an infinity algebra}
An $A_\infty$-algebra $V$ is called \emph{unital} if there exists an  element $\tau$ of degree $-1$ in $\Sigma V^*$ which could be extended to a basis $\tau, {\bf t}$ so that in this basis $m$ has a form
\[m=A({\bf t})\partial_\tau+\sum_i B({\bf t})\partial_{t_i}+\ad\tau-\tau^2\partial_\tau.\]
Note that in the dual basis of $V$ the element $\tau$ corresponds to the unit element.
Since $\tau^2=1/2[\tau,\tau]$ this definition also makes sense for $C_\infty$-algebras.
\begin{rem}
For a unital $A_\infty$-algebra $(T\Sigma V^*, m)$ the differential $m$ is always \emph{exact}, cf. for example, \cite{HL}, Lemma 6.9; this corresponds to the well-known fact that a (co)bar-construction of a unital associative algebra is contractible. A similar remark applies to unital $C_\infty$-algebras as well.
\end{rem}
 Suppose that $V$ is a unital augmented differential graded algebra. In other words $A$ admits a decomposition $V=k\cdot 1\oplus V_+$ where $V_+$ is  a differential graded algebra without a unit; alternatively one can say that $V$ is obtained from $V_+$ by adjoining a unit. Consider an $A_\infty$ minimal model of $V_+$; choosing a basis in $H^\bullet(V_+)$ and the corresponding dual basis $t_i$ in $H^\bullet(\Sigma V^*)$ we could assume that this minimal model has the form $(k\langle\langle{\bf t}\rangle\rangle, m_+)$ where $m_+$ is a derivation of $k\langle\langle{\bf t}\rangle\rangle$ of degree $-1$ with vanishing constant and linear terms. Then we have the following result.
\begin{prop}\label{aug1}
Under the above assumptions the derivation \\ $m:=m_++\ad \tau-\tau^2\partial_\tau$ of $k\langle\langle{\tau,\bf t}\rangle\rangle$ is a minimal unital $A_\infty$-model of $V$.
\end{prop}
\begin{proof}
Choose a basis  $\{x_i\}$ in $V_+$ and the corresponding dual basis $\{x^\prime_i\}$ in $\Sigma V^*$. Then the cobar-construction $\hat{T}\Sigma V_+^*$ of $V_+$ could be identified with the ring of noncommutative power series $k\langle\langle {\bf{x}}^\prime\rangle\rangle$. Denote by $m^\prime$ the differential in  $k\langle\langle {\bf{x}}^\prime\rangle\rangle$. 

The collection $[1, \{x_i\}]$ is a basis in $V$; consider $[\tau^\prime, \{x^\prime_i\}]$, the corresponding basis in $\Sigma V^*$. Then the cobar-construction  $\hat{T}\Sigma V^*$ of $V$ could be identified with the ring 
of noncommutative power series $k\langle\langle\tau^\prime, {\bf{x}}^\prime\rangle\rangle$ with the differential
$m^\prime-\tau^{\prime 2}\partial_\tau\prime+\ad\tau^\prime $.

Since $k\langle\langle{\bf t}\rangle\rangle$ endowed with the differential $m_+$ is a minimal model of $V_+$ there is a (continuous) map $f: k\langle\langle{\bf t}\rangle\rangle\rightarrow k\langle\langle{\bf x}^\prime\rangle\rangle$
which determines a  quasi-isomorphism $\Sigma H^\bullet(V^*_+)\rightarrow \Sigma V_+^*$
Then $f$ can be extended to the map of dga's 
\[\tilde{f}: k\langle\langle{\tau,\bf t}\rangle\rangle\rightarrow k\langle\langle{\tau^\prime,\bf x}^\prime\rangle\rangle\]
by letting $\tilde{f}(\tau)=\tau^\prime$.
It is clear that $\tilde{f}$ determines a quasi-isomorphism $\Sigma H^\bullet(V^*)\rightarrow \Sigma V^*$ as required. 
\end{proof}
We now formulate the analogue of this result in the $C_\infty$-context.  Let $V$ be a unital augmented \emph{commutative} differential graded algebra; $V=k\oplus V_+$. Consider a $C_\infty$ minimal model of $V_+$; choosing a basis in $H^\bullet(V_+)$ and the corresponding dual basis $t_i$ in $\Sigma V^*$ we could assume that this minimal model is the pair $(k\{\{{\bf t}\}\}, m_+)$ where $m_+$ is a derivation of $k\{\{{\bf t}\}\}$ of degree $-1$ with vanishing linear term. Note that a $C_\infty$ minimal model of a commutative differential graded algebra could be obtained from its $A_\infty$ minimal model by taking the primitive elements of the latter. Also note that $\tau^2\partial_\tau$ is a Lie derivation: $\tau^2\partial_\tau=1/2[\tau,\tau]\partial_\tau.$ With these remarks we have the following corollary.
\begin{cor}\label{aug2}
Under the above assumptions the derivation \\ $m:=m_++\ad \tau-1/2[\tau,\tau]\partial_\tau$ of $k\{\{{\tau,\bf t}\}\}$ is a $C_\infty$ minimal model of $V$.
\end{cor} \noproof                      
 Thus, our results provide minimal canonical unital models for augmented (associative or commutative) dga's.  Later on whenever we talk about minimal models of augmented (commutative) dga's we will mean these canonical models.

Next we consider morphisms. Recall that an $A_\infty$-morphism $\phi:k\langle\langle\tau,{\bf t}\rangle\rangle\rightarrow k\langle\langle\tau^\prime,{\bf t}^\prime\rangle\rangle$ between two unital $A_\infty$-algebras is called \emph{unital} if it has the form
\[\phi(\tau)=\tau^\prime+A(\bf t),\]
\[\phi(t_i)=B_i(\bf t).\] A $C_\infty$-morphism between two $C_\infty$-algebras is called \emph{unital} if the corresponding $A_\infty$-morphism is unital.

We have the following obvious result.
\begin{prop}\label{obvious}
Let $V=k\oplus V_+$ and $U=k\oplus U_+$ be augmented $A_\infty$-algebras, $k\langle\langle\tau,\bf t\rangle\rangle$ and  $k\langle\langle\tau^\prime,\bf t^\prime\rangle\rangle$ be their minimal unital models.  Then any $A_\infty$-morphism $V_+\rightarrow U_+$ determines a unital $A_\infty$-morphism
 \[k\langle\langle\tau^\prime,{\bf t}^\prime\rangle\rangle\rightarrow k\langle\langle\tau,\bf t\rangle\rangle.\] Similarly if $V=k\oplus V_+$ and $U=k\oplus U_+$ are augmented $C_\infty$-algebras and $k\{\{\tau,\bf t\}\}$ and  $k\{\{\tau^\prime,\bf t^\prime\}\}$ are their minimal unital models then any $C_\infty$-morphism $V_+\rightarrow U_+$ determines a unital $C_\infty$-morphism
 \[k\{\{\tau^\prime,{\bf t}^\prime\}\}\rightarrow k\{\{\tau,\bf t\}\}.\] 
\end{prop}
\noproof
In view of the above results it makes sense to define the adjunction of a unit to an arbitrary $A_\infty$ or $C_\infty$-algebra as follows.
\begin{prop}
Let $V$ be an $A_\infty$-algebra determined by the derivation $m:k\langle\langle {\bf x}\rangle\rangle\rightarrow k\langle\langle {\bf x}\rangle\rangle$. Then the derivation $\tilde{m}$ of the algebra $k\langle\langle {\bf x}, \tau\rangle\rangle$ determined by the formula
\[\tilde{m}=m+\ad\tau-\tau^2\partial_\tau\]
is a unital $A_\infty$-structure on $\tilde{V}:=k\oplus V$. We will call $\tilde{V}$ the $A_\infty$-algebra obtained from $V$ by adjoining a unit.
\end{prop}\noproof
It is clear that the above proposition-definition remains valid in the $C_\infty$-context and is independent of the choice of a basis. Moreover the correspondence $V\mapsto \tilde{V}$ is a functor from the category of $A_\infty$-algebras ($C_\infty$-algebras) to the category of unital $A_\infty$-algebras ($C_\infty$-algebras respectively).

\subsection{$C_\infty$-models for rational homotopy types}

Let $M$ be a nilpotent $CW$-complex of finite type and let $V:=A^\bullet(M)$ be the minimal Sullivan algebra of $M$.  
Recall that $V$ is a free graded commutative algebra with a decomposable differential $d$ (i.e. $d(V_+)\subset V_+\cdot V_+)$ that is multiplicatively quasi-isomorphic to the Sullivan-deRham algebra of $M$, cf. \cite{BG}. A Quillen model of $M$ can be identified with the space of primitives inside $T\Sigma V^*_+$, the cobar-construction of $V$; it inherits the differential  from $T\Sigma V^*_+$ and becomes a dgla.
\begin{defi}
A minimal Quillen model  of $M$ is a $C_\infty$ minimal model of the commutative dga $V_+=A^\bullet(M)_+$.
It will be defined by ${\mathcal L}(M)$.
\end{defi}  
\begin{rem}
Our notion of a minimal dgla differs slightly from that introduced in \cite{Nei}, \cite{BL} in that we consider \emph{completed} free Lie algebras with a decomposable differential. On the other hand if $M$ is simply-connected then ${\mathcal L}(M)$ has generators in strictly positive degrees, so the completion does not make any difference and the differential applied to each generator is always a finite sum of monomials.
Therefore in the simply-connected case our definition agrees with that of \cite{Nei} and \cite{BL}. In the nilpotent case the existence of a (conventional) minimal Quillen model is unknown, but our $C_\infty$-model always exists and provides a perfectly adequate substitute. However a real challenge would be to construct a $C_\infty$ minimal model encoding a \emph{nonnilpotent} rational homotopy type. In this connection note that, as shown by Neisendorfer \cite{Nei} a nonnilpotent dgla (e.g. a semisimple Lie algebra sitting in degree zero) in general does not admit a conventional minimal model.
\end{rem}
\begin{defi}
A contractible Quillen minimal model $\tilde{\mathcal L}(M)$ is a unital $C_\infty$ minimal model of $A^\bullet(M)$.
\end{defi}
\begin{rem}
Clearly a contractible Quillen minimal model $\tilde{\mathcal L}(M)$ is obtained from a nonunital $C_\infty$-algebra ${\mathcal L}(M)$ by the procedure of adjoining a unit as discussed in the previous subsection.
\end{rem}
From now on we will only consider the simply-connected case.  The next result shows that ${\mathcal L}(M)$ and $\tilde{\mathcal L}(M)$ faithfully record the rational homotopy type of $M$.
\begin{theorem}The following conditions are equivalent:\begin{enumerate}\item
Two simply-connected spaces $M$ and $N$ of finite type are rationally equivalent; \item
${\mathcal L}(M)$ and ${\mathcal L}(N)$ are isomorphic;\item
$\tilde{{\mathcal L}}(M)$ and $\tilde{{\mathcal L}}(N)$ are isomorphic through a unital $C_\infty$-isomorphism.\end{enumerate}
\end{theorem}  
\begin{proof}
The equivalence of (1) and (2) is well-known, cf. \cite{Nei} or \cite{BL}. The implication (2)$\Rightarrow$(3) is is Proposition \ref{obvious}. For (3)$\Rightarrow$(1) let \[\phi:\tilde{{\mathcal L}}(N)=k\{\{\tau^\prime,{\bf t}^\prime\}\}\rightarrow k\{\{\tau,{\bf t}\}\}=\tilde{\mathcal L}(M)\] be a unital $C_\infty$-isomorphism.
Then $\phi(\tau^\prime) =\tau+A(\bf t)$ but since $A(\bf t)$ has degree $\geq 0$ and $|\tau|=|\tau^\prime|=-1$ we conclude that ${A(\bf t)}=0$ and so $\phi$ restricts  to an isomorphism between
${\mathcal L}(M)$ and ${\mathcal L}(N)$. 
\end{proof}
\section{Cyclic cohomology of infinity-algebras.}
 \begin{defi}Let $V$ be an $A_\infty$-algebra and consider $C^\bullet_\lambda(V):=\Sigma (\hat{T}\Sigma V^*)_+/[,]$, the (suspension of the) quotient of the reduced completed tensor algebra by the space of all graded commutators. The derivation $m$ determines a differential on $C^\bullet_\lambda(V)$ making it into a complex. This complex will be called the cyclic complex of the $A_\infty$-algebra $V$ and its cohomology $HC^\bullet(V)$ the cyclic cohomology of $A$.\end{defi}
\begin{rem}
In more familiar terms the cyclic complex of $A$ is the complex of the form
\[V^*\rightarrow (V\otimes V)^*/[,]\rightarrow\ldots\rightarrow [V^*]^{\otimes n}/[,]\ldots \]
where the differential is determined by the $A_\infty$-structure $m$ and reduces to the familiar Connes complex, cf. \cite{Loday} if one disregards the higher products $m_i, i>2$ .
Note that the quotient $[V^*]^{\otimes n}/[,]$ is naturally identified with $[V^*]^{\otimes n}_{Z_n}$, the space of coinvariants with respect to the action of the cyclic group $Z_n$. The complex $C^\bullet_\lambda(V)$ is contravariantly functorial with respect to $A_\infty$-morphisms.
\end{rem}
\begin{prop}
A weak equivalence $V\rightarrow U$ between two $A_\infty$-algebras induces an isomorphism
$HC^\bullet(U)\rightarrow HC^\bullet(V)$.
\end{prop}
\begin{proof}
The complexes $C^\bullet_\lambda(V)$ and $C^\bullet_\lambda(U)$ have filtrations induced by wordlength in the tensor algebras  $(\hat{T}\Sigma V^*)_+$ and $(\hat{T}\Sigma U^*)_+$. Since the functor of $Z_n$-coinvariants is exact we conclude that the induced map on the $E_1$-terms of the corresponding spectral sequences is an isomorphism.
\end{proof}
Now consider the $A_\infty$-algebra $\tilde{V}$ obtained from $V$ by adjoining a unit. Then $V$ is an $A_\infty$-retract of $\tilde{V}$ which implies that $HC^\bullet(V)$ is a direct summand in $HC^\bullet(\tilde{V})$. More precisely, we have the following result.
\begin{prop}\label{isom}
There is a canonical isomorphism $HC^\bullet(\tilde{V})\cong HC^\bullet(V)\oplus HC^\bullet(k)$.
\end{prop}
\begin{proof}
This result is well-known in the case when $V$ is a dga. See \cite{Loday} for the proof in the ungraded case which is carried over almost verbatim to the dga case. Let $U$ be a dga which is $A_\infty$-equivalent to $V$. We have a commutative square of $A_\infty$-algebras whose horizontal maps
are weak equivalences and therefore induce isomorphisms in cyclic cohomology:
\[\xymatrix{V\ar[d]\ar[r]&U\ar[d]\\ \tilde{V}\ar[r]&\tilde{U}}\]
Then the desired isomorphism for $V$ follows from the corresponding result for $U$.
\end{proof}

\section{Poincar\'e duality spaces and symplectic infinity-algebras}
We start by recalling the notion of a symplectic (or cyclic) infinity-algebra. More details could be found in \cite{HL}.

Let $(\hat{T}\Sigma V^*,m)$ be an $A_\infty$-algebra. We assume that $V$ is finite dimensional over $k$ and that
$A$ possesses a nondegenerate graded symmetric scalar product $\langle,\rangle$ which we will refer to as the \emph{inner product}. Then $\Sigma V$ also acquires a scalar product which we will denote by the same symbol $\langle,\rangle$; namely:
\[\langle\Sigma a,\Sigma b\rangle:=(-1)^{|a|}\langle a,b\rangle.\]
It is easy to check that the product on $\Sigma V$ will be graded skew-symmetric, in other words it will determine a (linear) graded symplectic structure on $\Sigma V$.

We will consider the scalar product $\langle,\rangle$ as an element $\omega\in(\Sigma V^*)^{\otimes 2}$.  Consider the element $\tilde{m}:=m(\omega)\in \hat{T}\Sigma V^*$. Clearly $\tilde{m}=\tilde{m_1}+\tilde{m_2}+\ldots$ where $\tilde{m_i}$ has wordlength $i+1$. In other words, the tensors $\tilde{m_i}$ is obtained from $m_i$ by `raising an index' with the help of the form $\langle,\rangle$.

The space $T^n(\Sigma V^*)$ has an action of the symmetric group $S_n$ permuting the tensor factors. Note that each time a pair of elements $a,b$ in a monomial is permuted the result acquires the sign $(-1)^{|a||b|}$.
\begin{defi}
An $A_\infty$-algebra $(\hat{T}\Sigma V^*, m)$ with an inner product $\langle,\rangle$ is called \emph{symplectic} if $\tilde{m}$ is invariant with respect to all cyclic permutations of its tensor summands.
A $C_\infty$-algebra is called symplectic if it is so considered as an $A_\infty$-algebra.
\end{defi}
Given a basis $x_i$ in $\Sigma V^*$ the element $\omega\in T^2(\Sigma V^*)$ could be written as $\omega=\sum\omega^{ij}x_i\otimes x_j$. Consider the element \[[\omega]:=\sum\omega^{ij}[x_i,x_j]\in L^2(\Sigma V^*)\hookrightarrow T^2(\Sigma V^*).\] Clearly $[\omega]$ does not depend on the choice of a basis.  An easy calculation establishes the following result.
\begin{prop}
An $A_\infty$-algebra $(\hat{T}\Sigma V^*,m)$ (or $ C_\infty$-algebra $(\hat{L}\Sigma V^*,m)$) with an inner product $\omega$ is symplectic if and only if $m([\omega])=0$.
\end{prop}
\noproof
\begin{example}
Let $V$ be the cohomology algebra of $S^n$ with its usual Poincar\'e duality form.
It is easy to see that the corresponding symplectic $A_\infty$-algebra has the form $k\langle\langle \tau,y\rangle\rangle$ with $m=\ad\tau-\tau^2\partial_\tau, \omega=(-1)^{n}\tau\otimes y+y\otimes\tau$ and $[\omega]=2[y,\tau]$.
\end{example}   
Now suppose that $M$ is a simply-connected rational Poincar\'e duality space of dimension $n$. Consider a Quillen minimal model $\mathcal{L}(M)$ and the corresponding contractible Quillen model $\tilde{\mathcal L}(M)$. The $C_\infty$-algebra $\tilde{\mathcal L}(M)$ may not be symplectic, however it is so `up to homotopy'. In other words the quadratic part of its differential preserves the graded symplectic form on the space of its decomposables. This is just a reformulation of the invariance property of the Poincar\'e pairing $\langle ab, c\rangle =\langle a, bc\rangle$ on $H^\bullet(M)$. Then the main theorem of \cite{HL} states that $\tilde{\mathcal L}(M)$ is isomorphic to a \emph{symplectic} $S_\infty$-algebra, i.e. the whole differential, not just its quadratic part, preserves our symplectic form. Slightly modifying the proof in the cited reference one can show that the resulting symplectic $C_\infty$-algebra can be chosen to be unital. 
\begin{defi}
We will call a symplectic unital minimal $C_\infty$-algebra isomorphic to $\tilde{\mathcal L}(M)$  a \emph{canonical contractible model} of $M$.\end{defi}
 Let us denote a canonical contractible model of $M$ temporarily by $\tilde{\mathcal L}_c(M)$. Choose a basis in $H^\bullet(M)$ which includes $1\in H^0(M)$ and $[M]^*\in H^n(M)$, the dual to the fundamental cycle of $M$. We will denote the corresponding basis in $[\Sigma H^\bullet(M)]^*=\Sigma^{-1}H_{-\bullet}(M)$ by $\tau,\{x_i\},y$ where $\tau$ is the dual to $1$ and $y=[M]$.  Then the $C_\infty$-structure $m$ (the differential)  on   $\tilde{\mathcal L}_c(M)$ will have the following form:
\[m=\ad\tau-1/2[\tau,\tau]\partial_\tau+A({\bf x},y)\partial_\tau +\sum B_i({\bf x},y)\partial_{x_i}+C({\bf x},y)\partial_y.\] 
Since $m$ has degree $-1$ we conclude from dimensional considerations that $A({\bf x}, y)=0$ and that $B_i({\bf x},y)$ and $ C({\bf x},y)$ do not depend on $y$ so we can write $B_i({\bf x},y)=B_i({\bf x})$ and   $ C({\bf x},y)= C({\bf x})$. 

It follows that $\sum B_i({\bf x})\partial_{x_i}+C({\bf x})\partial_y$ determines a differential on ${k}\{\{{\bf x},y\}\}$ so we obtain a (nonunital) minimal $C_\infty$-algebra ${\mathcal L}_c(M)$. Clearly an isomorphism between $\tilde{\mathcal L}_c(M)$ and $\tilde{\mathcal L}(M)$ restricts to an isomorphism between ${\mathcal L}_c(M)$ and ${\mathcal L}(M)$ so ${\mathcal L}_c(M)$ could serve as a (minimal) Quillen model of $M$. We will call it a Stasheff model of $M$. From now on we will suppress the subscript $c$ for a Stasheff model and a canonical contractible model of $M$ since only those will be considered later on. 

Furthermore note that since $B_i(\bf x)$ does not depend on $y$ the derivation $\sum B_i({\bf x})\partial_{x_i}$ restricted to the Lie algebra $k\{\{\bf x\}\}$ has square zero and so determines a minimal $C_\infty$-algebra. This minimal $C_\infty$-algebra is a Quillen model of the $n-1$-skeleton $\dot{M}$ of $M$. It will be called a \emph{Stasheff model of $\dot{M}$}. 

Now consider the Poincar\'e duality form $\langle,\rangle$ on $H^\bullet(M)$. It is clear that it is the direct sum of its restrictions on the subspaces $H^0(M)\oplus H^n(M)$ and $\oplus_{0<i<n}H^i(M)$. Consequently the canonical element $[\omega]\in L^2(H_\bullet(M))$ could be represented as $[\omega]=[\bar{\omega}]+2[y,\tau]$ where $[\bar\omega]$ corresponds to the nondegenerate scalar product on $\oplus_{0<i<n}H^i(M)$. Note that $[\bar{\omega}]$ only involves the commutators of the $x_i$'s. The following result shows that the condition that $m([\omega])=0$ imposes further restrictions on $B_i(\bf x)$ and $C(\bf x)$.
\begin{theorem} Let $(k\{\{{\bf x},y\}\}, \ad\tau-1/2[\tau,\tau]\partial_\tau+\sum B_i({\bf x})\partial_{x_i}+C({\bf x})\partial_y)$ be a Stasheff model for a Poincar\'e duality space $M$. Then 
\begin{enumerate}\item $C(\bf x)$ is purely quadratic; in fact $C({\bf x})=1/2[\bar{\omega}]$.
\item The corresponding Stasheff model of $\dot{M}$ is symplectic;\\ in other words
$\sum B_i({\bf x)}\partial _{x_i}([\bar\omega])=0$.
       \end{enumerate}
\end{theorem}
\begin{proof}\
\begin{enumerate}\item
It is clear that $C({\bf x})=1/2[\bar{\omega}]+C^\prime(\bf x)$ where $C^\prime(\bf x)$ denote a sum of terms of order $>3$.
 Computing the part of $m(1/2[\bar{\omega}]-[\tau,y])$ containing the elements of bracket length $>3$ we see that it equal to $[\tau, C^\prime(\bf x)]+$ terms involving $x_i$'s only. It follows that $[\tau, C^\prime({\bf x})]=0$ and therefore $C^\prime({\bf x})=0$ as claimed.
\item  We have
\[0=2\cdot m\circ m(y)=2\cdot m(1/2[\omega])=\sum B_i({\bf x)}\partial _{x_i}([\bar\omega]).\]
 \end{enumerate}
\end{proof}
\begin{rem}
The last theorem was proved (in a different langauge) by Stasheff in \cite{Sta}. A correction of Stasheff's argument was later given in \cite{Aubry}, see also \cite{Umble}.
\end{rem}  
{\section {String bracket on the equivariant homology of the free loop space of $ \dot{M}$.}  }
\subsection{ Cyclic cohomology of symplectic $A_\infty$-algebras} Let $(\hat{T}\Sigma V^*, m)$ be a symplectic $A_\infty$-algebra and $[\omega]\in T^2\Sigma V^*$ be the canonical element corresponding to the invariant inner product on $V$. Consider the graded Lie algebra $\Der(\hat{T}\Sigma V^*)$ of all \emph{continuous} derivations of $\hat{T}\Sigma V^*$. The commutator with $m$ determines a differential $\Der(\hat{T}\Sigma V^*)$ making it into a (Hochschild) complex $C^\bullet(V)$. Denote by $SC^\bullet(V,V)$ the subcomplex of $C^\bullet(V,V)$ formed by \emph{symplectic} derivations, i.e. the derivations vanishing on $[\omega]$. 

Recall from \cite{HL} that $SC^\bullet(V,V)$ is  isomorphic to $\Sigma^{|\omega|-2}C_\lambda^\bullet(V)$, the cyclic Hochschild complex 
computing the cyclic cohomology of $(\hat{T}\Sigma V^*, m)$. 

Clearly the commutator of derivations determines the structure of a dgla on $SC^\bullet(V,V)$. Using the above isomorphism we obtain the bracket
\begin{equation}\label{bracket}HC^n(V)\otimes HC^l(V)\rightarrow HC^{n+l-|\omega|+1}(V).\end{equation}
\begin{example}
Let $V$ be a graded vector space together with a graded symmetric nondegenerate even or odd pairing $V\otimes V\rightarrow k$. Taking $m:\hat{T}\Sigma V^*\rightarrow \hat{T}\Sigma V^*$ to be the zero map we can view $V$ as a symplectic $A_\infty$-algebra. Then $H^\bullet_\lambda(V)$ could then be identified with the (completion of the) space of all cyclic words in $\Sigma V^*$. The bracket (\ref{bracket}) coincides with the one defined by Kontsevich \cite{Kon} in the context of his noncommutative symplectic geometry. 
\end{example}

\subsection { Equivariant homology of ${\mathbb L}M$ and cyclic cohomology}
Let $M$ be a simply-connected space of finite type. We denote by ${\mathbb L}M$ the space of unbased loops in $M$ and by $H^{S_1}_\bullet({\mathbb L}M):=ES^1\times_{S^1}{\mathbb L}M$ the $S^1$-equivariant homology of ${\mathbb L}M$.  There is a fibration  $ES^1\times_{S^1}{\mathbb L}M\rightarrow BS^1$ which determines a map of graded coalgebras \[ H^{S_1}_\bullet({\mathbb L}M)\rightarrow H_\bullet(BS^1)=({k}[u])^*.\]
Since this fibration has a section the above coalgebra map is split,  so we can identify $H_\bullet(BS^1)$ with its image in $H^{S_1}_\bullet({\mathbb L}M)$.
\begin{defi}
The \emph{reduced} $S^1$-equivariant homology $\bar{H}^{S_1}_\bullet({\mathbb L}M)$ of ${\mathbb L}M$ is defined as $\bar{H}^{S_1}_\bullet({\mathbb L}M):=H^{S_1}_\bullet({\mathbb L}M)/H_\bullet(BS^1)$.
\end{defi}
Recall from \cite{jones}, \cite{HL} that $H^{S_1}_\bullet({\mathbb L}M)$ could be expressed in terms of the cyclic cohomology of the cochain algebra of $M$:
\[ H^{S_1}_n({\mathbb L}M)\cong HC^{-n+1}(C^\bullet(M)).\]
The choice of a basepoint in $M$ gives $C^\bullet(M)$ an augmentation and so its cyclic cohomology contains a copy of the $HC^\bullet(k)=(k[u])^* $. We have therefore
\[HC^\bullet(C^\bullet(M))\cong HC^\bullet(k)\oplus \overline{HC}^\bullet(C^\bullet(M))\]
where  $\overline{HC}^\bullet(C^\bullet(M))$ is the \emph{reduced} cyclic cohomology of $C^\bullet(M)_+$ (the above isomorphism could be taken as a definition of the reduced cyclic cohomology in the augmented case).
It is clear that there is an isomorphism
\[ \bar{H}^{S_1}_n({\mathbb L}M)\cong \overline{HC}^{-n+1}(C^\bullet(M)).\]
We now assume that $M$ be a simply-connected rational Poincar\'e duality space of dimension $n$. Recall that we denoted by $\dot{M}$ the $n-1$-skeleton of $M$.  
\begin{theorem}
The equivariant homology of the loop space on $\dot{M}$ possesses the structure of a graded Lie algebra of degree $2-n$. The Lie bracket on $\bar{H}^{S_1}_\bullet({\mathbb L}\dot{M})$ will be referred to as the \emph{string bracket}. If $N$ is another Poincar\'e duality space homotopy equivalent to $M$ through an orientation-preserving homotopy equivalence then the corresponding graded Lie algebras are isomorphic.
\end{theorem}
\begin{proof}
Let $A^\bullet$ be the Sullivan minimal model of $\dot{M}$ and $A^\bullet_+$ be its space of indecomposable elements; clearly $A^\bullet$ is obtained from $A^\bullet_+$ by adjoining a unit. 
We have by Proposition \ref{isom} 
\begin{align*}\label{po} HC^\bullet(A^\bullet)&\cong HC^\bullet(A^\bullet_+)\oplus HC^\bullet(k)\\
&\cong \overline{HC}^\bullet(A^\bullet)\oplus HC^\bullet(k).\end{align*}
We see that $\overline{HC}^\bullet(A^\bullet)\cong \bar{H}^{S_1}_{-\bullet+1}({\mathbb L}\dot{M})$ is canonically identified with $HC^\bullet(A^\bullet_+)$. Recall that $A^\bullet_+$ has a minimal $C_\infty$-model that is symplectic -- a Stasheff model. Denote this model by $V$; we therefore have an isomorphism $HC^\bullet(A^\bullet_+)\cong HC^\bullet(V)$.
Using this isomorphism  the bracket (\ref{bracket}) could be transferred to $\bar{H}^{S_1}_\bullet({\mathbb L}\dot{M})$ as required. The homotopy invariance of the string bracket thus defined is evident.
\end{proof}
\begin{example}
Consider the wedge of $2N$ spheres of the form $X=\bigvee_{i=1}^{N}(S^{n_i}\vee S^{n-n_i})$ where $1<n_i<N-1$. Then we could build a Poincar\'e duality space $M$ out of $X$ by attaching an $n$-cell. We conclude that $X$ is homotopy equivalent to $\dot{M}$. Then the reduced equivariant cohomology of ${\mathbb L}\dot{M}$ can be identified with the cyclic cohomology of the corresponding zero-multiplication algebra which is isomorphic to the space of the cyclic words in $2N$ letters. The string bracket is Kontsevich's noncommutative Poisson bracket \cite{Kon}. 
\end{example}
\begin{rem}
Since ${H}^{S_1}_\bullet({\mathbb L}{M})$ could be identified (with an appropriate shift) with $HC^\bullet(C^\bullet(M))$ one can define a string bracket on  ${H}^{S_1}_\bullet({\mathbb L}{M})$ by taking cyclic cohomology of the (contractible) symplectic $A_\infty$-model of $C^\bullet(M)$. This was the approach of \cite{HL} and it is likely that the obtained string bracket agrees with that of Sullivan-Chas. Recall that this model has the form $k\langle\langle \tau, {\bf x}, y\rangle\rangle$. Since any symplectic derivation of  $k\langle\langle {\bf x} \rangle\rangle$ clearly extends to $k\langle\langle \tau, {\bf x}, y\rangle\rangle$ we conclude that the string brackets on $ {\mathbb L}{M}$ and ${\mathbb L}{\dot{M}}$ are compatible in the sense that the inclusion $\dot{M}\hookrightarrow M$ determines a map of corresponding graded Lie algebras. It would be interesting to give a geometric description of the string bracket on ${\mathbb L}{\dot{M}}$ along the lines of Sullivan-Chas.
 \end{rem} 
\begin{rem}
It appears that there are no corresponding analogues for the Chas-Sullivan loop product and loop bracket on the homology of the loop space of $\dot{M}$.
\end{rem}

\end{document}